\documentclass[11pt,a4paper]{article}
\usepackage[T1]{fontenc}
\IfFileExists{lmodern.sty}{\usepackage{lmodern}}{}
\usepackage{amsmath,amssymb,amsthm}
\usepackage[a4paper,margin=29mm]{geometry}
\usepackage[hidelinks]{hyperref}
\newtheorem{theorem}{Theorem}[section]
\newtheorem{lemma}[theorem]{Lemma}
\newtheorem{corollary}[theorem]{Corollary}
\theoremstyle{remark}
\newtheorem{remark}[theorem]{Remark}

\newcommand{\R}{\mathbb R}
\newcommand{\one}{\mathbf 1}
\newcommand{\norm}[1]{\lVert #1\rVert}
\newcommand{\ip}[2]{\langle #1,#2\rangle}
\newcommand{\cV}{\mathcal V}
\newcommand{\cA}{\mathcal A}

\title{Small subset sums in Euclidean space:\\ the sharp constant and its extremal families}
\author{Dami\'an Pinasco\\[3pt]\small Departamento de Matem\'atica y Estad\'istica, Universidad Torcuato Di Tella\\\small CONICET, Buenos Aires, Argentina\\\small \texttt{dpinasco@utdt.edu}}
\date{}

\begin{document}
\maketitle

\begingroup
\renewcommand{\thefootnote}{}
\footnotetext{%
}
\endgroup

\begin{abstract}
Let $v_1,\ldots,v_m$ be vectors in the Euclidean unit ball
of $\R^n$ with sum zero. For every $k\in\{1,\ldots,m\}$,
one can select exactly $k$ vectors whose sum has squared
norm at most $L_n=\lfloor(n+1)^2/4\rfloor/n$.
We characterize all families attaining this bound in
dimensions $n\ge3$: they consist of equally many copies
of the vertices of a regular simplex inscribed in the
unit sphere, with the possibility of one additional
zero vector in even dimensions. We specify the admissible
cardinalities. We also obtain a sharp mean-square rounding
bound with prescribed marginals, for a law that can be
chosen independently of the vectors. The sharp deterministic
bound follows from the affine reduction of Ambrus,
B\'ar\'any and Grinberg and the rounding lemma of Tinarrage.
Following Tinarrage's variational strategy, we give a direct
geometric proof of the latter inequality and determine its
equality cases for intermediate cardinalities.
\end{abstract}

\section{Introduction}

A basic problem in finite-dimensional geometry is to select vectors whose sum is small. Ambrus, B\'ar\'any and Grinberg \cite{abg} studied the following version. For the unit ball $B$ of a norm on $\R^n$, let $T(B)$ be the least number such that every finite family $v_1,\ldots,v_m\in B$ with $\sum_iv_i=0$ and every integer $0\le k\le m$ admit a set $A\subseteq\{1,\ldots,m\}$ with $|A|=k$ and $\norm{\sum_{i\in A}v_i}\le T(B)$. Families are indexed, so repetitions are allowed. They proved $T(B)\le\lceil n/2\rceil$ for every norm, with equality for suitable norms, and that the Euclidean constant is of order $\sqrt n$ \cite[Theorem 4]{abg}.

If $u_1,\ldots,u_{n+1}$ are the vertices of a centered regular simplex inscribed in the unit sphere, then $\ip{u_i}{u_j}=-1/n$ for $i\ne j$, and every sum of $k$ distinct vertices has squared norm $k(n+1-k)/n$. Its maximum over $k$ gives the lower bound in \cite{abg}. The affine reduction of \cite{abg} and the sharp rounding lemma of Tinarrage \cite[Lemma 7]{tinarrage} yield the matching upper bound. We determine the full equality cases.

\begin{theorem}\label{th:main}
\textup{(i)} For the Euclidean unit ball $B_2^n$,
\[
 T(B_2^n)^2=L_n:=\frac{\lfloor(n+1)^2/4\rfloor}{n}.
\]
Equivalently, $T(B_2^n)=\sqrt{n+2}/2$ when $n$ is even and $T(B_2^n)=(n+1)/(2\sqrt n)$ when $n$ is odd.

\medskip
\noindent\textup{(ii)} Let $n\ge3$, $p=n+1$ and $r=\lfloor p/2\rfloor$, and let $v_1,\ldots,v_m\in B_2^n$ have sum zero. For an integer $0\le k\le m$, $\min_{|A|=k}\norm{\sum_{i\in A}v_i}^2=L_n$ if and only if, up to an orthogonal transformation and a permutation, one of the following holds.
\begin{enumerate}
\item The family consists of $c\ge1$ copies of each vertex of a centered regular simplex inscribed in the unit sphere, and $k\equiv r$ or $k\equiv p-r\pmod p$.
\item The dimension $n$ is even, the family consists of $c\ge1$ copies of each vertex of such a simplex together with one zero vector, and $k\equiv r+1\pmod p$.
\end{enumerate}
\end{theorem}

The first part is a consequence of \cite{abg,tinarrage}. The ideas of the proofs are described below.

We also consider a probabilistic form of the rounding problem, with prescribed marginals.

For integers $p\ge1$ and $0\le q\le p$ put
\[
 \Delta(p,q)=\Bigl\{x\in[0,1]^p:\ \sum_i x_i=q\Bigr\},\qquad
 \cV(p,q)=\Delta(p,q)\cap\{0,1\}^p .
\]
The second set is the vertex set of the first: if two coordinates are fractional they can be moved in opposite directions, and a point with exactly one fractional coordinate cannot have integer coordinate sum.

For $p\ge2$ put
\[
 R(p,q)=\frac{q(p-q)}{p-1}.
\]

\begin{theorem}\label{th:law}
Let $p\ge2$ and $q\in\{0,\ldots,p\}$. For every $x\in\Delta(p,q)$ there are weights $\mu_y\ge0$, $y\in\cV(p,q)$, with $\sum_y\mu_y=1$ and $\sum_y\mu_yy=x$, such that
\[
 \sum_y\mu_y\Bigl\|\sum_i(y_i-x_i)u_i\Bigr\|^2\le R(p,q)
\]
for every family $u_1,\ldots,u_p$ of vectors of norm at most one in a real Hilbert space. If $2\le q\le p-2$ and $x\ne(q/p)\one$, the right-hand side can be replaced by a strictly smaller number depending only on $x$, $p$ and $q$.
\end{theorem}

The same law works for every family of vectors. Its construction uses convex envelopes and the classical minimax principle.

\subsection*{Idea of the proof}
\emph{The reduction.} By the affine reduction of \cite{abg} (Lemma~\ref{lem:affine}), after fixing some of the vectors one is left with at most $n+1$ vectors $u_1,\ldots,u_p$, coefficients $x_i\in(0,1)$ with integer sum $q$, and the task of choosing exactly $q$ of these vectors so that the error $\sum_i(y_i-x_i)u_i$ is small. The key inequality (Theorem~\ref{th:round}, due to Tinarrage) says that some choice has squared error at most $R(p,q)$. Geometrically: every point of the convex hull of the sums of $q$ of the vectors $u_i$ lies within distance $\sqrt{R(p,q)}$ of one of these sums. Maximizing $R(p,q)$ over $p\le n+1$ gives $L_n$, and the regular simplex shows that nothing is lost.

\emph{The rounding inequality.} Consider a worst configuration, that is, coefficients and vectors for which the best choice is as bad as possible, and call a choice active if it attains this worst value $t$. At such a configuration no small simultaneous change of the coefficients and of the vectors can make all active choices worse at once. A separation argument turns this into a probability distribution on the active choices that is in equilibrium with respect to both kinds of changes. Equilibrium for the coefficients says that the distribution has mean $x$: it is a way of rounding $x$ at random. Equilibrium for the vectors says that each $u_i$ points in a direction determined by the others, and linear algebra then forces the covariance of the random choice into a rigid form: there is a probability vector $\theta$ with $\sum_i\theta_iu_i=0$ such that the $i$th coordinate of the random choice has variance $x_i(1-x_i)$ equal to $t$ times the variance $\theta_i(1-\theta_i)$ of a Bernoulli variable of parameter $\theta_i$. Since unit vectors in equilibrium cannot carry a weight larger than one half, every $\theta_i\le1/2$. Taking square roots, each coordinate contributes a plane vector whose length is proportional to $1-2\theta_i$; these lengths add up to a fixed total, whereas the vectors themselves add up to a vector determined by $p$, $q$ and $t$. The triangle inequality in the plane then gives $t\le R(p,q)$. Equality forces $\theta_i=1/p$ and $x_i=q/p$ for every $i$, and an averaging argument identifies the regular simplex.

\emph{The extremal families.} If an original family attains $L_n$, all the inequalities above are equalities, so the reduction produces a regular simplex formed by $n+1$ members of the family, with equal coefficients. The sums of a fixed number of its vertices span a polytope containing a ball about the origin of a definite radius; adding any further vector of the family to these sums must not decrease the error, and this forces every further vector to be a vertex of the same simplex or, in even dimensions, the zero vector. The zero-sum condition and the prescribed cardinality then determine the multiplicities and the possible exceptional zero.

\emph{The common law.} For a fixed family, the least mean-square error of a random rounding with mean $x$ is a convex-envelope quantity; at a point where it is largest it equals the distance to the nearest admissible choice, a phenomenon going back to Wegmann \cite{wegmann}. The rounding inequality therefore bounds it by $R(p,q)$, and a minimax argument produces a single law that works for all families simultaneously.

The rounding inequality used in the proof is the barycentric inequality of Tinarrage \cite[Lemma 7]{tinarrage}, established in connection with a sharp approximate Carath\'eodory theorem and iterated Delaunay refinement. Following the overall strategy of \cite{tinarrage}, we give a geometric proof with a different derivation of the contact identity, using direct variations of the vectors and finite-dimensional separation, and a final step based on the triangle inequality in the plane. The equality analysis provides the starting point for the classification.

The equilibrium conditions in the proof are reminiscent of John's contact identities \cite{john,ball}. Related balance and isotropy conditions arise in extremal problems for signed sums of unit vectors \cite{ambrusnietert,pinasco}. In our setting, the zero-sum condition identifies the sum of $k$ selected vectors with half a signed sum subject to $\sum_i\varepsilon_i=2k-m$. Thus Theorem~\ref{th:main}(i) is also a sharp cardinality-constrained vector balancing result.

\section{Preliminaries}

\subsection{Affine reduction}

We begin with a reduction due to Ambrus, B\'ar\'any and Grinberg \cite{abg}. Instead of selecting a subset directly, we allow fractional coefficients while preserving both their sum and the vector sum. Passing to an extreme point of the resulting polytope leaves at most $n+1$ fractional coefficients. Thus the number of vectors still to be rounded is controlled by the dimension, independently of the size of the original family.

\begin{lemma}\label{lem:affine}
Let $v_1,\ldots,v_m\in\R^n$ and $c\in\R^n$, and let $z^0\in[0,1]^m$ satisfy $\sum_iz^0_i=k$ and $\sum_iz^0_iv_i=c$. There are disjoint sets $A_0,J\subseteq\{1,\ldots,m\}$ and numbers $0<x_i<1$ $(i\in J)$ such that
\[
 \sum_{i\in A_0}v_i+\sum_{i\in J}x_iv_i=c,\qquad \sum_{i\in J}x_i=k-|A_0|,\qquad |J|\le n+1 .
\]
\end{lemma}

\begin{proof}
The polytope $\{z\in[0,1]^m:\sum_iz_i=k,\ \sum_iz_iv_i=c\}$ contains $z^0$. At one of its extreme points let $A_0$ index the coordinates equal to one and $J$ the fractional ones. The vectors $(v_i,1)\in\R^{n+1}$, $i\in J$, are linearly independent, since a dependence would allow a small perturbation of the fractional coordinates in both directions. Hence $|J|\le n+1$.
\end{proof}

\subsection{A separation principle}

The following elementary form of finite-dimensional convex separation will be used to obtain a common system of weights for the variations of an extremal configuration. We include its proof for completeness.

\begin{lemma}\label{lem:sep}
Let $K\subset\R^N$ be compact and convex, and suppose that no point of $K$ has all coordinates positive. Then there is $\pi\in\R^N$ with $\pi_i\ge0$, $\sum_i\pi_i=1$, and $\ip{\pi}{k}\le0$ for every $k\in K$.
\end{lemma}

\begin{proof}
Let $E=K-\R^N_{\ge0}=\{k-o:\ k\in K,\ o\ge0\}$. It is convex, and closed because $K$ is compact. Let $e=(1,\ldots,1)$ and $\varepsilon>0$. The point $\varepsilon e$ is not in $E$: if $\varepsilon e=k-o$, then $k=\varepsilon e+o$ would have all coordinates positive. Let $z\in E$ be a point nearest to $\varepsilon e$ and $\pi^\varepsilon=(\varepsilon e-z)/\norm{\varepsilon e-z}$. For $w\in E$ and $0<s\le1$ the point $z+s(w-z)$ lies in $E$, so $\norm{\varepsilon e-z-s(w-z)}^2\ge\norm{\varepsilon e-z}^2$; expanding and letting $s\to0$ gives $\ip{\pi^\varepsilon}{w-z}\le0$. With $w=z-o$, $o\ge0$, this gives $\pi^\varepsilon\ge0$. With $w=k\in K$ it gives
\[
 \ip{\pi^\varepsilon}{k}\le\ip{\pi^\varepsilon}{z}<\ip{\pi^\varepsilon}{\varepsilon e}\le\varepsilon\sqrt N .
\]
A limit point $\pi'$ of $\pi^\varepsilon$ as $\varepsilon\to0$ satisfies $\pi'\ge0$, $\norm{\pi'}=1$ and $\ip{\pi'}{k}\le0$ on $K$. Normalize $\pi=\pi'/\sum_i\pi'_i$.
\end{proof}

\subsection{A minimax principle}

The common rounding law requires a choice of weights that works simultaneously for all families of vectors. For this purpose we use the following special case of the classical minimax theorem \cite{sion}, which follows directly from the preceding separation principle.

\begin{lemma}\label{lem:minimax}
Let $\mathcal P$ be a polytope with vertices $\mu^1,\ldots,\mu^M$, let $\mathcal K$ be a compact convex subset of a finite-dimensional space, and let $\Phi$ be affine in each variable. Put $v=\max_{G\in\mathcal K}\min_{\mu\in\mathcal P}\Phi(\mu,G)$. Then there is $\mu\in\mathcal P$ with $\Phi(\mu,G)\le v$ for all $G\in\mathcal K$.
\end{lemma}

\begin{proof}
For each $G$, the minimum over $\mathcal P$ of the affine function $\Phi(\cdot,G)$ is attained at a vertex, so $\min_j\Phi(\mu^j,G)\le v$; this minimum is continuous in $G$, so $v$ is attained. Fix $\varepsilon>0$. The set $\{(\Phi(\mu^j,G)-v-\varepsilon)_{j=1}^M:\ G\in\mathcal K\}$ is compact and convex, being an affine image of $\mathcal K$, and none of its points has all coordinates positive. Lemma~\ref{lem:sep} gives $\pi^\varepsilon\ge0$ with $\sum_j\pi^\varepsilon_j=1$ and $\Phi(\mu_\varepsilon,G)\le v+\varepsilon$ for all $G$, where $\mu_\varepsilon=\sum_j\pi^\varepsilon_j\mu^j\in\mathcal P$. A limit point of $\mu_\varepsilon$ as $\varepsilon\to0$ has the required property.
\end{proof}

\section{The rounding inequality}\label{sec:round}

\begin{theorem}\label{th:round}
Let $u_1,\ldots,u_p$ be vectors of norm at most one in a real Hilbert space, $p\ge2$, and $x\in\Delta(p,q)$ with $q$ an integer. Then there exists $y\in\cV(p,q)$ with
\begin{equation}\label{eq:round}
 \Bigl\|\sum_i(y_i-x_i)u_i\Bigr\|^2\le R(p,q).
\end{equation}
The bound is attained when $x=(q/p)\one$ and the $u_i$ are the vertices of a centered regular simplex inscribed in the unit sphere. If $2\le q\le p-2$, these are the only cases in which the minimum over $y$ of the left-hand side equals $R(p,q)$.
\end{theorem}

The inequality is due to Tinarrage \cite[Lemma 7]{tinarrage}. The following proof uses variations of the vectors and the coefficients and includes the equality analysis needed for Theorem~\ref{th:main}(ii).

\begin{remark}\label{rem:tin}
The rounding inequality is equivalent to Tinarrage's barycentric inequality \cite[Lemma 7]{tinarrage}.
For $1\le q\le p-1$, since $\sum_i(y_i-x_i)=0$, the error vector in \eqref{eq:round} does not change if all the $u_i$ are translated by the same vector. The $u_i$ lie in their affine hull $E$, of dimension at most $p-1$, and $E$ meets the unit ball in a ball of radius at most one. Tinarrage's lemma, in the version for points in a ball of radius $\delta\le1$ (his proof allows affinely dependent points), therefore gives \eqref{eq:round} with $d=p-1$, $r=q$ and barycentric coordinates $\lambda_i=x_i/q$. Conversely, his lemma is the special case of \eqref{eq:round} in which the points lie in $\R^{p-1}$.
\end{remark}

The proof follows the outline given in the introduction: Section~\ref{ss:prelim} reduces to interior coefficients and affinely independent vectors, Section~\ref{ss:weights} produces the equilibrium distribution, Section~\ref{ss:contact} derives the contact identity, Section~\ref{ss:triangle} applies the triangle inequality, and Section~\ref{ss:equality} treats equality.

\subsection{Initial reductions}\label{ss:prelim}
The assertion is trivial for $q=0$ and $q=p$. For $q=1$ put $a=\sum_ix_iu_i$; then
\[
 \sum_ix_i\norm{u_i-a}^2=\sum_ix_i\norm{u_i}^2-\norm a^2\le1,
\]
so one of the errors $u_i-a$ has squared norm at most one. Replacing $x$ and $y$ by $\one-x$ and $\one-y$ handles $q=p-1$. For the regular simplex and $x=(q/p)\one$, every error has squared norm $q-q(q-1)/(p-1)=R(p,q)$.

Let $2\le q\le p-2$; thus $p\ge4$. We argue by induction on $p$, simultaneously for all $q$. A family of $p$ vectors can be represented isometrically in $\R^p$, so we may maximize the continuous function
\[
 F(x,u)=\min_{y\in\cV(p,q)}\Bigl\|\sum_i(y_i-x_i)u_i\Bigr\|^2
\]
over the compact set $\Delta(p,q)\times(B_2^p)^p$. Let $(x,u)$ be a maximizer and $t=F(x,u)$. The regular simplex gives $t\ge R(p,q)>1$. We shall prove $t\le R(p,q)$ and analyze the case of equality; since every pair at which the minimum equals $R(p,q)$ is a maximizer, this proves the theorem.

\emph{Every coordinate of $x$ lies strictly between $0$ and $1$.} If $x_i=0$, restrict the choices to $y_i=0$; the remaining coordinates form a point of $\Delta(p-1,q)$, and the induction hypothesis bounds $F(x,u)$ by $R(p-1,q)$. If $x_i=1$, we obtain $R(p-1,q-1)$ in the same way. Both are strictly smaller than $R(p,q)$, because
\[
 R(p,q)-R(p-1,q)=\frac{q(q-1)}{(p-1)(p-2)},\qquad
 R(p,q)-R(p-1,q-1)=\frac{(p-q)(p-q-1)}{(p-1)(p-2)} .
\]

\emph{The vectors $u_i$ are affinely independent.} Otherwise there are numbers $h_i$, not all zero, with $\sum_ih_i=0$ and $\sum_ih_iu_i=0$. Moving $x$ along $x+sh$ until a coordinate reaches $0$ or $1$ changes no error vector, and gives a maximizer with a coordinate in $\{0,1\}$, which is impossible.

\subsection{Active weights}\label{ss:weights}
Let $\cA$ be the set of $y\in\cV(p,q)$ at which the minimum defining $t$ is attained, and write
\[
 a_y=\sum_i(y_i-x_i)u_i\qquad(y\in\cA).
\]
For $z\in\Delta(p,q)$ and $w=(w_1,\ldots,w_p)\in(B_2^p)^p$ move from $(x,u)$ towards $(z,w)$ along the segment. The derivative at the starting point of $\norm{a_y}^2$ is
\begin{equation}\label{eq:deriv}
 d_y(z,w)=2\Bigl\langle a_y,\ \sum_i(y_i-x_i)(w_i-u_i)-\sum_i(z_i-x_i)u_i\Bigr\rangle .
\end{equation}
The vectors $(d_y(z,w))_{y\in\cA}$ form a compact convex subset of $\R^{\cA}$, since \eqref{eq:deriv} is affine in $(z,w)$. None of them has all coordinates positive: otherwise all active errors would increase along a short common step, while the finitely many inactive ones stay above $t$, contradicting maximality. By Lemma~\ref{lem:sep} there are numbers $\pi_y\ge0$ with $\sum_{y\in\cA}\pi_y=1$ and
\begin{equation}\label{eq:sep}
 \sum_{y\in\cA}\pi_y\,d_y(z,w)\le0\qquad\text{for all }(z,w).
\end{equation}
The essential point is that the same weights control the variations of the coefficients and of the vectors.

\emph{The weights have barycenter $x$.} Take $w=u$. Since $x$ is interior, $z=x-s(\bar y-x)$ is admissible for small $s>0$, where $\bar y=\sum_y\pi_yy$. As $\sum_y\pi_ya_y=\sum_i(\bar y_i-x_i)u_i$, inequality \eqref{eq:sep} gives $2s\norm{\sum_i(\bar y_i-x_i)u_i}^2\le0$. Since $\sum_i(\bar y_i-x_i)=0$, affine independence yields
\begin{equation}\label{eq:mean}
 \sum_{y\in\cA}\pi_yy_i=x_i\qquad(1\le i\le p).
\end{equation}

\subsection{The contact identity}\label{ss:contact}
Put
\[
 c_{ij}=\sum_{y\in\cA}\pi_y(y_i-x_i)(y_j-x_j).
\]
The matrix $(c_{ij})$ is symmetric, its row sums vanish because $\sum_i(y_i-x_i)=0$, and by \eqref{eq:mean} its diagonal entries are $x_i(1-x_i)$.

Take $z=x$ in \eqref{eq:sep} and vary only the $i$th vector. With $b_i=\sum_jc_{ij}u_j$ we get $\ip{b_i}{w_i-u_i}\le0$ for every $w_i$ in the unit ball. A linear functional attains its maximum on the unit ball only in its own direction, so there are numbers $\lambda_i\ge0$ with
\begin{equation}\label{eq:normal}
 \sum_jc_{ij}u_j=\lambda_iu_i,\qquad \lambda_i(1-\norm{u_i}^2)=0
\end{equation}
(take $\lambda_i=\norm{b_i}$). Taking the scalar product with $u_i$ and summing over $i$,
\begin{equation}\label{eq:sumlambda}
 \sum_i\lambda_i=\sum_{i,j}c_{ij}\ip{u_i}{u_j}=\sum_{y\in\cA}\pi_y\norm{a_y}^2=t .
\end{equation}

Let $S$ be the symmetric matrix with entries $s_{ij}=\lambda_i\delta_{ij}-c_{ij}$. By \eqref{eq:normal}, each column of $S$ is the coefficient vector of a linear relation among the $u_i$. As the $u_i$ are affinely independent, these relations form a space of dimension at most one. On the other hand $S\one=\lambda$, because the rows of $(c_{ij})$ sum to zero, and $\lambda\ne0$ by \eqref{eq:sumlambda}. Hence every column of $S$ is a multiple of $\lambda$; write $s_{ij}=\lambda_ir_j$. Symmetry forces $r$ to be proportional to $\lambda$, and $S\one=\lambda$ with $\sum_i\lambda_i=t$ gives $s_{ij}=\lambda_i\lambda_j/t$. With $\theta_i=\lambda_i/t$ we have proved
\begin{equation}\label{eq:contact}
 \sum_i\theta_i=1,\qquad \sum_i\theta_iu_i=0,\qquad c_{ij}=t\,(\theta_i\delta_{ij}-\theta_i\theta_j);
\end{equation}
the second identity holds because $\lambda=S\one$ is itself a linear relation. On the diagonal,
\begin{equation}\label{eq:scalar}
 x_i(1-x_i)=t\,\theta_i(1-\theta_i).
\end{equation}
Since $0<x_i<1$, every $\theta_i$ is positive, so $\norm{u_i}=1$ by \eqref{eq:normal}. Finally, by \eqref{eq:contact},
\[
 \theta_i=\norm{\theta_iu_i}=\Bigl\|\sum_{j\ne i}\theta_ju_j\Bigr\|\le1-\theta_i,
\]
that is, $\theta_i\le1/2$.

\subsection{A triangle inequality}\label{ss:triangle}
Since $\theta_i\le1/2$, identity \eqref{eq:scalar} can be written as
\[
 \sqrt t\,(1-2\theta_i)=\sqrt{t-1+(2x_i-1)^2}=\bigl\|(\sqrt{t-1},\,2x_i-1)\bigr\| .
\]
Summing over $i$ and applying the triangle inequality in $\R^2$, with $\sum_i\theta_i=1$ and $\sum_ix_i=q$,
\[
 \sqrt t\,(p-2)\ \ge\ \bigl\|(p\sqrt{t-1},\,2q-p)\bigr\|=\sqrt{p^2(t-1)+(2q-p)^2}.
\]
Squaring, $t(p-2)^2-p^2(t-1)-(2q-p)^2=4(p-1)\bigl(R(p,q)-t\bigr)\ge0$. Thus $t\le R(p,q)$, which completes the induction.

\subsection{Equality}\label{ss:equality}
Suppose $t=R(p,q)$. Then equality holds in the triangle inequality, so the vectors $(\sqrt{t-1},2x_i-1)$ are positive multiples of each other. Their first coordinates coincide and are positive, hence all $x_i=q/p$. By \eqref{eq:scalar} and $\theta_i\le1/2$ all $\theta_i$ are equal to $1/p$, and \eqref{eq:contact} shows that the $u_i$ are unit vectors with $\sum_iu_i=0$.

Choose a $q$-subset uniformly at random. An index is included with probability $q/p$ and two given indices with probability $q(q-1)/(p(p-1))$; using $\sum_iu_i=0$ one finds that the mean squared norm of the selected sum is $R(p,q)$. As the minimum is $R(p,q)$, every sum of $q$ of the vectors has squared norm $R(p,q)$. For distinct $i,j$ and a set $E$ of size $q-1$ avoiding them, comparing the sums over $E\cup\{i\}$ and $E\cup\{j\}$ gives
\[
 \sum_{\ell\in E}\bigl(\ip{u_\ell}{u_i}-\ip{u_\ell}{u_j}\bigr)=0 .
\]
Since $1\le q-1\le p-3$, exchanging one element of $E$ shows that all summands, as $\ell$ ranges outside $\{i,j\}$, coincide; their sum over $E$ vanishes, so they vanish. Hence all off-diagonal scalar products are equal, and $\sum_iu_i=0$ fixes their common value at $-1/(p-1)$: the $u_i$ form a centered regular simplex. This proves Theorem~\ref{th:round}.

\begin{remark}\label{rem:compare}
Tinarrage's proof of \cite[Lemma 7]{tinarrage} also maximizes jointly over the configuration and the point, applies Gordan's alternative to the active barycenters \cite[Lemma 6]{tinarrage}, and identifies the second moment of the active weights; identity \eqref{eq:contact} corresponds to \cite[(30)]{tinarrage}, and \eqref{eq:scalar} to \cite[(44)]{tinarrage}, his vector $\tau_G$ being our $\theta$. The arguments differ in two places. First, we vary the vectors themselves rather than their Gram matrix. This avoids the description of first-order Gram deformations on the sphere \cite[Lemma 5]{tinarrage} and the positivity argument for $\tau_G$; the identity $S=\lambda\lambda^\top/t$ follows from the fact that each column of $S$ is a linear relation, and the bound $\theta_i\le1/2$ from the equilibrium $\sum_i\theta_iu_i=0$. Second, starting from \eqref{eq:scalar} we conclude with a triangle inequality in the plane instead of the analysis of \cite[Subcases 3c(i)--(ii)]{tinarrage}. This last step also gives the equality cases. Tinarrage's reduction to unit vectors \cite[Lemma 8]{tinarrage} corresponds to our use of \eqref{eq:normal}.
\end{remark}

\begin{remark}
The structure of the proof of Theorem~\ref{th:round} resembles the contact conditions in John's theorem \cite{john,ball}: separation provides weights, and the variations of the vectors give an equilibrium. By \eqref{eq:contact}, for every vector $a$,
\[
 \sum_{y\in\cA}\pi_y\ip{a_y}{a}^2=t\sum_i\theta_i\ip{u_i}{a}^2 ,
\]
an identity between two quadratic forms that holds before the extremal configuration is identified. At the regular extremal, if $L=\operatorname{span}\{u_1,\ldots,u_p\}$ and $P_L$ is the orthogonal projection onto $L$, both sides equal $t\norm{P_La}^2/(p-1)$. Thus the second moment is isotropic on $L$.
\end{remark}

\section{The sharp bound and its consequences}

\begin{proof}[Proof of Theorem~\ref{th:main}(i)]
Let $v_1,\ldots,v_m\in B_2^n$ have sum zero and $0\le k\le m$. Apply Lemma~\ref{lem:affine} with $c=0$ and $z^0=(k/m)\one$, and put $p=|J|$, $q=k-|A_0|$. If $p=0$, the sum over $A_0$ vanishes; $p=1$ is impossible because $q$ is an integer. Otherwise Theorem~\ref{th:round}, applied to $(v_i)_{i\in J}$ and $x=(x_i)_{i\in J}$, gives $B\subseteq J$ with $|B|=q$ and
\[
 \Bigl\|\sum_{i\in A_0\cup B}v_i\Bigr\|^2=\Bigl\|\sum_{i\in B}v_i-\sum_{i\in J}x_iv_i\Bigr\|^2\le R(p,q).
\]
For fixed $q$, $R(p+1,q)-R(p,q)=q(q-1)/(p(p-1))\ge0$; as $p\le n+1$, $R(p,q)\le R(n+1,q)\le L_n$. The regular simplex gives the lower bound.
\end{proof}

Let $T(B_2^n,k)$ denote the analogous constant for a fixed cardinality $k$, the supremum being over families with at least $k$ members. The preceding argument gives the following consequences.

\begin{corollary}\label{cor:k}
If $1\le k\le\lfloor(n+1)/2\rfloor$, then $T(B_2^n,k)^2=k(n+1-k)/n$.
\end{corollary}

For $k=2$ this also follows from Jung's theorem \cite{jung}. The same argument gives the sharp Euclidean form of the corrected \cite[Theorem 7]{abg}, which allows an arbitrary target in the convex hull with coefficients at most $1/k$.

\begin{corollary}\label{cor:target}
Let $w_1,\ldots,w_m\in B_2^n$, $1\le k\le m$, and $z=\sum_i\lambda_iw_i$ with $\lambda_i\ge0$, $\sum_i\lambda_i=1$ and $\lambda_i\le1/k$. Then there is $A$ with $|A|=k$ and $\norm{\sum_{i\in A}w_i-kz}^2\le L_n$. The constant is sharp.
\end{corollary}

\begin{proof}[Proof of Corollary~\ref{cor:k}]
In the proof above $q\le k\le\lfloor(n+1)/2\rfloor$, and $R(n+1,\cdot)$ is increasing on $\{0,\ldots,\lfloor(n+1)/2\rfloor\}$. The regular simplex attains $k(n+1-k)/n$.
\end{proof}

\begin{proof}[Proof of Corollary~\ref{cor:target}]
Apply Lemma~\ref{lem:affine} with $c=kz$ and $z^0=k\lambda\in[0,1]^m$. If $J$ is empty, the sum over $A_0$ is already $kz$. Otherwise $p=|J|\ge2$, and the same argument gives $A=A_0\cup B$ with $\norm{\sum_{i\in A}w_i-kz}^2\le R(p,q)\le L_n$. Sharpness follows from the regular simplex with $z=0$ and $k=\lfloor(n+1)/2\rfloor$.
\end{proof}

For a fixed zero-sum family the minimum for cardinality $k$ equals the minimum for $m-k$, by taking complements. This symmetry alone does not imply that the largest minimum occurs at the central cardinality.

\section{Extremal families}

\begin{proof}[Proof of Theorem~\ref{th:main}(ii)]
Use an extreme point as in the proof of Theorem~\ref{th:main}(i). If the original minimum equals $L_n$, all inequalities in that proof are equalities. As $L_n>1$, the fractional cardinality $q$ and $p'-q$ are at least two, where $p'=|J|$. Since $R(\cdot,q)$ is strictly increasing for $q\ge2$, we get $p'=n+1=p$ and $q\in\{r,p-r\}$. By the equality case of Theorem~\ref{th:round}, the fractional vectors, denoted $u_1,\ldots,u_p$, form a centered regular simplex inscribed in the unit sphere and all their coefficients equal $q/p$. Consequently $\sum_{i\in A_0}v_i=0$.

Let $v$ be a vector of the family whose coefficient at the extreme point is zero. For $E\subseteq\{1,\ldots,p\}$ with $|E|=s:=q-1$ put $u_E=\sum_{i\in E}u_i$. Selecting $v$, the vectors indexed by $E$ and those indexed by $A_0$ gives $k$ vectors, so
\begin{equation}\label{eq:outside}
 \norm{u_E+v}^2\ge R(p,q)\qquad(|E|=s).
\end{equation}
The convex hull $P_s$ of the points $u_E$ is
\[
 P_s=\Bigl\{z:\ -\frac{s}{p-1}\le\ip{z}{u_i}\le\frac{p-s}{p-1}\ \ (1\le i\le p)\Bigr\}.
\]
Indeed, every $z$ can be written uniquely as $z=\sum_ib_iu_i$ with $\sum_ib_i=s$, and then $\ip{z}{u_i}=(pb_i-s)/(p-1)$, so the inequalities say exactly that $0\le b_i\le1$.

Put $h=s/(p-1)$. Since $q$ is central, $s\le(p-1)/2$, so $h\le1/2$ and $P_s$ contains the ball of radius $h$ about the origin. Hence $\min_{|E|=s}\ip{v}{u_E}\le-h\norm v$. As $\norm{u_E}^2=R(p,s)$ and $R(p,s+1)-R(p,s)=1-2h$, inequality \eqref{eq:outside} gives
\[
 (\norm v-1)(\norm v+1-2h)\ge0 .
\]
Since $\norm v\le1$, either $\norm v=1$, or $v=0$ and $h=1/2$. If $\norm v=1$, equality holds in the estimate for the minimum, so $-hv$ lies on the boundary of $P_s$. The upper constraints are at distance $(p-s)/(p-1)>h$ from the origin, hence a lower constraint is active at $-hv$; this means $\ip{v}{u_i}=1$ for some $i$, that is, $v=u_i$.

For a vector indexed by $A_0$, apply the same argument to the complementary extreme point $\one-z$ and the cardinality $m-k$: by the zero sum the minimum is again $L_n$, the fractional vectors are the same and the fractional cardinality is $p-q$. Thus every vector of the family is a vertex of the simplex or zero.

The only linear relation among the vertices is $\sum_iu_i=0$, so the zero sum forces a common multiplicity $c\ge1$. If $j$ nonzero vectors are selected, with $t_i$ copies of $u_i$, the squared norm of their sum is $\frac p{p-1}\sum_i(t_i-j/p)^2$. Under $0\le t_i\le c$ and $\sum_it_i=j$ this is minimized when the $t_i$ differ by at most one (move a unit from a larger to a smaller count if they differ by at least two), and such counts are feasible for $0\le j\le cp$. If $a$ is the residue of $j$ modulo $p$, the minimum is $R(p,a)$.

If the family contains $z_0$ zero vectors, the possible values of $j$ are the integers from $\max(0,k-z_0)$ to $\min(k,cp)$. Equality with $L_n>0$ excludes $j=0$ and $j=cp$, so this range is $[k-z_0,k]$, and all its residues must be central. For even $p$ there is one central residue and for odd $p$ two consecutive ones. Hence either $z_0=0$ and $k$ is central modulo $p$, or $p$ is odd, $z_0=1$ and $k\equiv r+1\pmod p$. Conversely, these families attain $L_n$ by the same computation.
\end{proof}

The extremal families also determine the constant for infinitely many large cardinalities.

\begin{corollary}\label{cor:large}
Let $n\ge3$, $p=n+1$ and $r=\lfloor p/2\rfloor$. If $k\equiv r$ or $k\equiv p-r\pmod p$, then $T(B_2^n,k)^2=L_n$. For every other $k\ge1$, no family attains the value $L_n$ for the cardinality $k$.
\end{corollary}

\begin{proof}
The upper bound is Theorem~\ref{th:main}(i). For the admissible residues, the families of Theorem~\ref{th:main}(ii), with $c$ large enough that the family has at least $k$ members, attain $L_n$. The last assertion is the `only if' part of Theorem~\ref{th:main}(ii).
\end{proof}

For $n=1,2$ one has $L_n=1$, and many other families are extremal.

\section{A common rounding law}

The quantity minimized below is the squared analogue of the effective standard deviation of Cassels \cite{cassels}; its relation to covering radii goes back to Wegmann \cite{wegmann}, see also \cite[Theorem 2.15]{fmmz}. Prescribed-marginal sampling itself is classical. In particular, Tinarrage \cite[Lemma 4]{tinarrage} uses the variance comparison for pivotal sampling of Chauvet and Ruiz-Gazen \cite{chauvet} to obtain a mean-square estimate. Here the sharp deterministic bound on each face and minimax give a common law with the sharper bound $R(p,q)$, which is strictly smaller than $\min(q,p-q)$ for $1<q<p-1$.

\begin{proof}
Fix a family $u$ and write $Uz=\sum_iz_iu_i$. For $z\in\Delta(p,q)$ let
\[
 g(z)=\min\Bigl\{\sum_y\mu_y\norm{Uy}^2:\ \mu_y\ge0,\ \sum_y\mu_y=1,\ \sum_y\mu_yy=z\Bigr\},\qquad h(z)=g(z)-\norm{Uz}^2 .
\]
The graph of $g$ is the lower boundary of the convex hull of the finitely many points $(y,\norm{Uy}^2)$, so $g$ is convex, piecewise affine and continuous; $h(z)$ is the least mean-square error of a law with barycenter $z$.

Let $z_0$ maximize $h$ and let $B$ be the smallest face of $\Delta(p,q)$ containing $z_0$. Every law with barycenter $z_0$ is supported on the vertices of $B$. If $B$ is a vertex, $h(z_0)=0$. Otherwise $z_0$ is relatively interior to $B$, and there is an affine function $\ell\le g$ on $B$ with $\ell(z_0)=g(z_0)$. The function $\ell(z)-\norm{Uz}^2\le h(z)$ attains its maximum over $B$ at $z_0$, so its derivative along the affine hull of $B$ vanishes there, and for every vertex $y$ of $B$
\[
 \ell(z_0)+2\ip{Uz_0}{U(y-z_0)}=\ell(y)\le g(y)=\norm{Uy}^2 .
\]
Rearranging, $h(z_0)\le\min_y\norm{U(y-z_0)}^2$ over the vertices $y$ of $B$. Theorem~\ref{th:round}, applied in the free coordinates of $B$, bounds this minimum by $R(p,q)$, since fixing coordinates never increases the bound (Section~\ref{sec:round}). If $2\le q\le p-2$, the bound on a proper face is strictly smaller than $R(p,q)$, so a maximizer with $h(z_0)=R(p,q)$ is interior and, by the equality statement of Theorem~\ref{th:round}, equals $(q/p)\one$.

To choose one law for all families, let $\mathcal K$ be the set of symmetric positive semidefinite $p\times p$ matrices with diagonal entries at most one; it is compact and convex, and consists of the Gram matrices of families of $p$ vectors in the unit ball of $\R^p$. For fixed $x$ the laws with barycenter $x$ form a polytope $\mathcal P_x$, and
\[
 \Phi(\mu,G)=\sum_y\mu_y\sum_{i,j}G_{ij}(y_i-x_i)(y_j-x_j)
\]
is bilinear. For each $G\in\mathcal K$, $\min_{\mu\in\mathcal P_x}\Phi(\mu,G)=h(x)\le R(p,q)$ by the first part. The classical minimax theorem \cite{sion} (in the elementary form of Lemma~\ref{lem:minimax}) gives $\mu\in\mathcal P_x$ with $\Phi(\mu,G)\le\max_{G'\in\mathcal K}\min_{\mu'\in\mathcal P_x}\Phi(\mu',G')\le R(p,q)$ for all $G$. In the strict case the function $G\mapsto\min_{\mu'}\Phi(\mu',G)$ is continuous and everywhere smaller than $R(p,q)$ on the compact set $\mathcal K$, hence its maximum is smaller than $R(p,q)$.
\end{proof}

The bound is sharp: take uniform marginals and the vertices of a centered regular simplex. Every admissible selection has squared error $R(p,q)$, so no choice of law can improve the bound uniformly.

\begin{remark}
The preceding theorem also yields a diagonal majorant for the covariance matrix.
For the law of Theorem~\ref{th:law} put $C_{ij}=\sum_y\mu_y(y_i-x_i)(y_j-x_j)$. There are numbers $d_i\ge0$ with $\sum_id_i\le R(p,q)$ and $\sum_{i,j}C_{ij}a_ia_j\le\sum_id_ia_i^2$ for all $a\in\R^p$. Indeed, maximize $\sum_{i,j}C_{ij}\ip{u_i}{u_j}$ over families in the unit ball of $\R^{p+1}$. As in \eqref{eq:normal}, a maximizer satisfies $\sum_jC_{ij}u_j=d_iu_i$ with $d_i\ge0$ and $d_i(1-\norm{u_i}^2)=0$, so $\sum_id_i$ equals the maximum, which is at most $R(p,q)$. A family of $p$ vectors spans at most $p$ dimensions, so there is a unit vector $e$ orthogonal to all $u_i$. Replacing $u_i$ by $\sqrt{1-s^2a_i^2}\,u_i+sa_ie$ keeps the vectors in the unit ball, and the coefficient of $s^2$ in the change of the objective is $\sum_{i,j}C_{ij}a_ia_j-\sum_id_ia_i^2$, which must therefore be nonpositive.
\end{remark}

\section{Larger cardinalities}

For $k\le\lfloor(n+1)/2\rfloor$ and for the residues of Corollary~\ref{cor:large}, the constant $T(B_2^n,k)$ is determined. For the remaining cardinalities and $n\ge3$, the value $L_n$ is not attained by any family (Corollary~\ref{cor:large}), although we do not know whether $T(B_2^n,k)^2$ is strictly smaller than $L_n$; moreover, the extremal families, if they exist, are not built from a single regular simplex. The following example in $\R^4$, with $k=4$, was found by a numerical search carried out with the assistance of an AI system and then identified in exact form; its verification is elementary.

In $\R^4$, with $\one=(1,1,1,1)$, put
\[
 a=-\tfrac12\,\one,\qquad b_i=\tfrac{\sqrt{13}}{4}\,e_i+\tfrac{5-\sqrt{13}}{16}\,\one\quad(1\le i\le4).
\]
Then $\norm a=\norm{b_i}=1$, $\ip{b_i}{b_j}=3/16$ for $i\ne j$, $\ip{a}{b_i}=-5/8$, and $\sum_ib_i=\frac54\one$. Consider the family consisting of three copies of $a$, the vector $-a/2$ and $b_1,\ldots,b_4$; its sum is zero. A selection of four members consists of $\alpha\le3$ copies of $a$, $\eta\in\{0,1\}$ copies of $-a/2$ and $\beta$ of the vectors $b_i$, with $\alpha+\eta+\beta=4$. With $c=\alpha-\eta/2$, the sum $S$ of the selection satisfies
\[
 \norm S^2=c^2-\tfrac54\,c\beta+\beta+\tfrac3{16}\,\beta(\beta-1).
\]
The eight possible triples $(\alpha,\eta,\beta)$ give $\norm S^2=11/8$ for $(1,0,3)$, $(2,0,2)$, $(1,1,2)$, $(2,1,1)$, and $\norm S^2=25/4$ for $(0,0,4)$, $(0,1,3)$, $(3,0,1)$, $(3,1,0)$. Hence
\[
 T(B_2^4,4)^2\ \ge\ \frac{11}{8},
\]
whereas copies of the vertices of a single regular simplex, possibly together with zero vectors, give at most $1$ for $k=4$, and the upper bound is $L_4=3/2$. It would be interesting to determine $T(B_2^n,k)$ for these cardinalities.

\section*{Use of AI tools}
During the preparation of this manuscript, the author used large language models (LLMs) for exploratory discussion, bibliographic assistance, consistency checks, numerical verification of explicit formulas, and language editing. The author has reviewed and verified the mathematical arguments, calculations, references and final text, and assumes full responsibility for the contents of the manuscript.

\end{document}